\input ./style/arxiv-vmsta.cfg
\documentclass[numbers,compress,v1.0.1]{vmsta}
\usepackage{vtexbibtags}
\usepackage{authorquery}

\volume{6}
\issue{1}
\pubyear{2019}
\firstpage{133}
\lastpage{144}
\aid{VMSTA118}
\doi{10.15559/18-VMSTA118}
\articletype{research-article}

\startlocaldefs

\allowdisplaybreaks

\newtheorem{theorem}{Theorem}[section]
\theoremstyle{definition}
\newtheorem{example}{Example}[section]
\endlocaldefs

\begin{document}

\begin{frontmatter}
\pretitle{Research Article}

\title{Ruin probability for the bi-seasonal discrete time risk model with dependent claims}


\author{\inits{O.}\fnms{Olga}~\snm{Navickien\.{e}}\ead[label=e1]{olga.navickiene@mif.vu.lt}}
\author{\inits{J.}\fnms{Jonas}~\snm{Sprindys}\ead[label=e2]{jonas.sprindys@mif.vu.lt}} %
\author{\inits{J.}\fnms{Jonas}~\snm{\v{S}iaulys}\thanksref{cor1}\ead[label=e3]{jonas.siaulys@mif.vu.lt}}
\thankstext[type=corresp,id=cor1]{Corresponding author.}

\address{Institute of Mathematics, \institution{Vilnius University},
Naugarduko 24, Vilnius LT-03225, \cny{Lithuania}}



\markboth{O.~Navickien\.{e} et al.}{Ruin probability for the bi-seasonal discrete time risk model with dependent claims}

\begin{abstract}
The discrete time risk model with two seasons and dependent claims is
considered. An algorithm is created for computing the values of the
ultimate ruin probability. Theoretical results are illustrated with
numerical examples.
\end{abstract}

\received{\sday{19} \smonth{7} \syear{2018}}
\revised{\sday{22} \smonth{9} \syear{2018}}
\accepted{\sday{22} \smonth{9} \syear{2018}}
\publishedonline{\sday{01} \smonth{10} \syear{2018}}
\end{frontmatter}

\vspace*{6pt}
\section{Introduction}%
\label{}
In this paper, we consider the bi-seasonal discrete time risk model\index{discrete time risk model} with
dependent claims.

\textit{We say that the insurer's surplus $W_{u}$ varies according to
the bi-seasonal risk model with dependent claims if}
\begin{equation*}
W_{u}(n)=u+n-\sum\limits
_{i=1}^{n}Z_{i}
\end{equation*}
\textit{for all} $n\in \mathbb{N}_{0}=\{0,1,2,\ldots \}$ \textit{and the following
assumptions hold}:
\begin{itemize}
\item
 \textit{the initial insurer's surplus is} $u\in \mathbb{N} _{0}$,
\item
\textit{there exists a random vector} $(X,Y)$ \textit{such that}
$(Z_{2k-1},Z_{2k})\,\mathop{=}\limits ^{d}\,(X,Y)$, $k\in \mathbb{N}$,
\item
\textit{the random vectors} $(Z_{2k-1},Z_{2k})$, $k\in
\mathbb{N}$, \textit{are independent},
\item
\textit{the generating random vector} $(X,Y)$ \textit{has the
distribution defined by the table below}, \textit{where} $h_{i,j} = \mathbb{P}$ $(X=i,
Y=j)$, $i,j \in \mathbb{N}_{0}$:
\begin{table}[h]
%
%
\begin{tabular}{cccccc}
\hline
$X\backslash Y$ & $0$ & $1$ & $2$ & $3$ & $\ldots $ \\
\hline
$0$ & $h_{0,0}$ & $h_{0,1}$ & $h_{0,2}$ & $h_{0,3}$ & $\ldots $ \\
$1$ & $h_{1,0}$ & $h_{1,1}$ & $h_{1,2}$ & $h_{1,3}$ & $\ldots $ \\
$2$ & $h_{2,0}$ & $h_{2,1}$ & $h_{2,2}$ & $h_{2,3}$ & $\ldots $ \\
$\ldots $ & $\ldots $ & $\ldots $ & $\ldots $ & $\ldots $ & $\ldots $ \\
\hline
\end{tabular}
%
%
\end{table}
\end{itemize}

If $X$ and $Y$ are independent random variables, then the model reduces
to the one considered in \cite{ds}. If, in addition, $X$ and
$Y$ are identically distributed, then the bi-seasonal discrete time risk
model\index{discrete time risk model} with dependent claims becomes the classical discrete time risk
model.\index{discrete time risk model}

The time of ruin and the ruin probability\index{ruin probability} are the main extremal
characteristics of insurance risk models. The time of ruin is defined
by the equality
\begin{equation*}
T_{u}= %
\begin{cases}
\min \{n\geqslant 1: W_{u}(n)\leqslant 0\},
\\
\infty , \ \textnormal{if}\ W_{u}(n)>0\ \textnormal{for\ all }\ n
\in \mathbb{N}.
\end{cases} %
\end{equation*}

The ultimate ruin probability,\index{ultimate ruin probability} or simply ruin probability,\index{ruin probability} is defined
by the following equality:
\begin{equation*}
\psi (u)=\mathbb{P}(T_{u}<\infty ).
\end{equation*}

In the case of the classical discrete time risk model,\index{discrete time risk model} recursive
procedures for calculating exact values of $\psi (u)$ are well known.
These procedures and related information can be found in
\cite{dg-1999,d-1994,d-2005,dw-1991,g-1988,gs-1998,llg-2009,s-1989,w-1993} among others.

The recursive calculation of $\psi (u)$ is relatively simple in the
classical discrete time risk model\index{discrete time risk model} because of the explicit formula for
$\psi (0)$. If the consecutive claim amounts $Z_{1}, Z_{2}, \ldots $ are
no longer identically distributed or independent, then the classical
discrete time risk model\index{discrete time risk model} becomes the inhomogeneous discrete time risk
model.\index{discrete time risk model} For all such models, the algorithms for finding values of the
ruin probabilities\index{ruin probability} are much more complicated. Several results related
to the calculation of the ruin probabilities for inhomogeneous renewal
risk models can be found in \cite{abp-2015,bl-2016,bs-2017,bs-2012,bbs-2010,ccglm-2013,cdnp-2016,ds,gks-2015,gks-2015a,hzyj-2017,rvz-2015,rvz-2015a,rvz-2017}
and~\cite{zfly-2017}.

The aim of this paper is to derive an algorithm for computing the values
of the ultimate ruin probability\index{ultimate ruin probability} in the bi-seasonal discrete time risk
model\index{discrete time risk model} with dependent claims. Theoretical results are illustrated with
numerical examples.

The rest of the paper is organized as follows. In Section~\ref{m}, we
present our main results. In Sections~\ref{p} and~\ref{pp}, the proofs
of the main results are given. Finally, in Section~\ref{e} we present
some examples, which show the applicability of our results.

\section{Main results}%
\label{m}
Let us introduce some notation used in our results. By
\begin{equation*}
x_{k}=\mathbb{P}(X=k), \qquad  y_{k}=\mathbb{P}(Y=k),\qquad
s_{k}=\mathbb{P}(S=k), \quad  k\in \mathbb{N}_{0},
\end{equation*}
we denote the marginal distributions\index{marginal distribution} of the random variables $X$,
$Y$ and their sum $S=X+Y$, respectively. The distribution functions of
these random variables are denoted by $F_{X}$, $F_{Y}$ and $F_{S}$, i.e.
\begin{gather*}
F_{X}(u)=\mathbb{P}(X\leqslant u)=\sum\limits
_{k=0}^{\lfloor u
\rfloor }x_{k}, \qquad  F_{Y}(u)=
\mathbb{P}(Y\leqslant u)=\sum\limits
_{k=0}^{\lfloor u\rfloor
}y_{k},
\\
F_{S}(u)=\mathbb{P}(S\leqslant u)=\sum
\limits
_{k=0}^{
\lfloor u\rfloor }s_{k}
\end{gather*}
for all $u\geq 0$. The notation $\overline{F}$ is used for the tail of
an arbitrary distribution function~$F$, i.e. $\overline{F}(u)=1-F(u)$
for all $u\in \mathbb{R}$.

Furthermore, the survival probability\index{survival probability} is denoted by $\varphi (u) = 1 -
\psi (u)$ for all $u \in \mathbb{N}_{0}$. It should be noted that our
main results are formulated in terms of the survival probability.\index{survival probability}

\begin{theorem}
\label{th1}
Let the bi-seasonal discrete time risk model\index{discrete time risk model} be generated by the random
vector $(X,Y)$,\index{random vector} where $X$ and $Y$ are nonnegative and integer-valued
random variables such that $\mathbb{E}X+\mathbb{E}Y<2$. In this case
%
\begin{equation}
\label{a0} \lim \limits
_{u\rightarrow \infty }\varphi (u)=1.
\end{equation}
\begin{itemize}
\item
 If $s_{0} = h_{0,0}>0$, then
%
\begin{gather}
\label{a1} \varphi (0) = (2-\mathbb{E}S) \lim \limits
_{n\rightarrow \infty }\ \frac{b
_{n+1}-b_{n}}{a_{n}-a_{n+1}},
\\
%
\label{a2} \varphi (u) = a_{u} \varphi (0) + b_{u} (2-
\mathbb{E}S), \quad  u\in \mathbb{N},
\end{gather}
where $a_{n}$ and $b_{n}$ are two sequences of real numbers defined
recursively by the equalities:
\begin{align*}
& a_{1}=-\frac{1}{y_{0}},\qquad  a_{n}=\frac{1}{s_{0}}
\Biggl(a_{n-2} - \sum\limits
_{i=1}^{n-1}s_{i}a_{n-i}+a_{1}
h_{n-1,0} \Biggr), \quad n\in \{2,3, \ldots \};
\\
& b_{1}=\frac{1}{y_{0}},\qquad  b_{n}= \frac{1}{s_{0}}
\Biggl(b_{n-2}-\sum\limits
_{i=1}^{n-1}s_{i}b_{n-i}+b_{1}
h_{n-1,0} \Biggr),\quad  n\in \{2,3, \ldots \}.
\end{align*}

\item If $s_{0} = 0$ with $x_{0} \neq 0$ and $y_{0} = 0$,
then
\begin{eqnarray*}
&& \varphi (0) = 2-\mathbb{E}S,
\\
&&\varphi (u) = \frac{1}{s_{1}} \Biggl(\varphi (u-1) - \sum
\limits
_{k=2} ^{u} s_{k} \varphi (u-k+1) \Biggr),\quad  u \in
\mathbb{N}.
\end{eqnarray*}

\item If $s_{0} = 0$ with $x_{0}=0$ and $y_{0} \neq 0$,
then
\begin{align*}
& \varphi (0) = 0,
\\
& \varphi (1) = \frac{1}{y_{0}}(2-\mathbb{E}S),
\\
& \varphi (u) = \frac{1}{s_{1}} \Biggl(\varphi (u-1) - \sum
\limits
_{k=2} ^{u} s_{k} \varphi (u-k+1) +
h_{u,0}\,\varphi (1) \Biggr),\quad  u\in \{ 2,3, \ldots \}.
\end{align*}
\end{itemize}
\end{theorem}

\begin{theorem}
\label{th2}
Let the bi-seasonal discrete time risk model\index{discrete time risk model} be generated by the random
vector $(X,Y)$,\index{random vector} where $X$ and $Y$ are nonnegative and integer-valued
random variables such that the net profit condition is not satisfied,
i.e. $\mathbb{E}X+\mathbb{E}Y \geqslant 2$.

If $\mathbb{E}X+\mathbb{E}Y > 2$, then $\varphi (u) = 0$ for all
$u \in \mathbb{N}_{0}$.

If $\mathbb{E}X+\mathbb{E}Y = 2$, then we have the following possible
subcases:
\begin{itemize}
\item
$\varphi (u) = 0$, $u \in \mathbb{N}_{0}$, {if} $s_{2} = h_{0,2} + h_{1,1} + h_{2,0} < 1$;
\item
$\varphi (0) = 0$, $\varphi (u) = 1$, $u \in \mathbb{N}$, {if}
$s_{2} = 1$ {and} $h_{2,0} = 0$;
\item
$\varphi (0) = \varphi (1) = 0$, $\varphi (u) = 1$, $u \in \{ 2,3,
\ldots \}$, {if} $s_{2} = 1$ {and} $h_{2,0} > 0$.
\end{itemize}
\end{theorem}

\section{Proof of Theorem~\ref{th1}}%
\label{p}
The proof is greatly influenced by the proofs given in
\cite{ds}. Therefore, many details that can be found there are omitted.

At the beginning of the proof consider the general case with
$\mathbb{E}S \geqslant 0$. By the total probability formula, we get the
following basic recursive formula for all $u \in \mathbb{N}_{0}$:
%
\begin{align}
\label{p1} \varphi (u) & = \sum\limits
_{k=0}^{u+1}
s_{k} \varphi (u+2-k) - h_{u+1,0} \varphi (1)
\nonumber
\\
& = \sum\limits
_{k=0}^{u+1} s_{u+1-k} \varphi (k+1) -
h_{u+1,0} \varphi (1).
\end{align}

The obtained equality implies that
\begin{eqnarray*}
\sum\limits
_{l=0}^{u} \varphi (l) = \sum
\limits
_{l=0}^{u} \sum\limits
_{k=0}^{l+1}
s_{l+1-k} \varphi (k+1) - \varphi (1) \sum\limits
_{l=0}^{u} h_{l+1,0}, \quad  u \in \mathbb{N}_{0}.
\end{eqnarray*}

By rearranging the terms we obtain
\begin{align*}
\sum\limits
_{k=0}^{u+2} \varphi (k)
\overline{F}_{S}(u+2-k) &= \varphi (u+1) + \varphi (u+2)
\\
&\quad  -  \varphi (1) \sum\limits
_{l=0}^{u+1}
h_{l,0} - \varphi (0) F_{S}(u+2).
\end{align*}

Passing to the limit as $u\rightarrow \infty $ in the last equality and
applying arguments similar to those in \cite{ds} we get
%
\begin{equation}
\label{p4} (2-\mathbb{E}S) \varphi (\infty ) = y_{0} \varphi (1)
+ \varphi (0).
\end{equation}

Now let us restrict to the case $\mathbb{E}S < 2$. Equality
\eqref{a0} is proved using the strong law of large numbers, and the
proof is identical to the proof of the first part of Theorem~2.3 in
\cite{ds}. As a result we get
%
\begin{equation}
\label{p5} 2-\mathbb{E}S = y_{0} \varphi (1) + \varphi (0).
\end{equation}

Suppose now that $s_{0} = h_{0,0} \neq 0$. Then \eqref{a2} can be
derived by induction with induction basis obtained from \eqref{p5}.
Equality \eqref{a1} can be derived in a way similar to that in
\cite{ds} with only the difference that the coefficients $a_{n}$ used
in the proof are different.

It remains to consider the case where $s_{0} = h_{0,0} = 0$. Since
$\mathbb{E}S < 2$, it follows that $s_{1} \neq 0$. Two subcases can be
considered separately: $x_{0} \neq 0$ and $y_{0} = 0$, or $x_{0} = 0$
and $y_{0} \neq 0$.

In the subcase where $x_{0} \neq 0$ and $y_{0} = 0$, we get the formula
for $\varphi (0)$ from \eqref{p5}. The formula for $\varphi (u)$,
$ u \in \mathbb{N}$, follows from \eqref{p1} because
\begin{equation*}
0 = y_{0} = \sum\limits
_{k=0}^{\infty }h_{k,0}
\end{equation*}
in the considered case.

If $x_{0} = 0$ and $y_{0} \neq 0$, then we get $\varphi (0)=0$ from
\eqref{p1}. Then the formula for $\varphi (1)$ follows from
\eqref{p5}, and the formula for $\varphi (u)$ in the case $ u \in \{ 2,3,
\ldots \}$ can be derived from \eqref{p1}.

Theorem~\ref{th1} is proved.

\section{Proof of Theorem~\ref{th2}}%
\label{pp}
Let us consider the cases $\mathbb{E}S>2$ and $\mathbb{E}S=2$
separately. The case $\mathbb{E}S>2$ can be proved using the same
arguments as in \cite{ds}.

In the case $\mathbb{E}S=2$, we can easily see from \eqref{p4} that
%
\begin{equation}
\label{p6} y_{0} \varphi (1) + \varphi (0) = 0.
\end{equation}
Therefore, $\varphi (0) = 0$. To calculate $\varphi (u), u \in
\mathbb{N}$, the subcases $s_{2}<1$ and $s_{2}=1$ can be considered
separately.

Consider the subcase $s_{2}<1$ first. We can prove that $\varphi (u) =
0,\ u \in \mathbb{N}$, in a way similar to that in \cite{ds} using
the fact that $\varphi (1) h_{l,0} = 0$ for $ l \in \mathbb{N}_{0}$,
which follows immediately from equality \eqref{p6}.

Now let us consider the subcase $s_{2} = h_{0,2} + h_{1,1} + h_{2,0} =
1$. There are the following possible cases:
\begin{itemize}
\item
If $h_{2,0} > 0$, then from the main recursive formula \eqref{p1} we get
$\varphi (1) = 0$.
\item
If $h_{2,0} = 0$, then obviously $W_{1}(n) \geqslant 1, n \in
\mathbb{N}$, and therefore, $\varphi (1) = 1$.
\end{itemize}

For $u \in \{2,3,\ldots \}$, it is easy to show that $W_{u}(n) \geqslant
1$ for $ n \in \mathbb{N}$, and therefore, $\varphi (u) = 1$ for such
$u$.

Theorem~\ref{th2} is proved.

\section{Numerical examples}%
\label{e}
In this section, four numerical examples for the calculation of the ruin
probability $\psi (u), u\in \mathbb{N}_{0}$,\index{ruin probability} are given. The first case
deals with the bivariate Poisson distribution,\index{bivariate Poisson distribution} and the next three cases
deal with a Clayton copula.\index{Clayton copula} The use of copulas is beneficial since it
gives the possibility of modeling marginal distributions\index{marginal distribution} and dependence
between them separately. Furthermore, while the bivariate Poisson
distribution\index{bivariate Poisson distribution} allows to model only positive dependence between marginals,\index{marginal}
a Clayton copula\index{Clayton copula} enables to model negative dependence\index{negative dependence} as well.

The numerical simulation procedure goes as follows. First, we can
calculate sufficiently many terms of the sequences $a_{u}$ and
$b_{u}$ from Theorem~\ref{th1}. Next, we can approximate $\psi (0)$ by
\begin{equation*}
\psi _{N}(0) = 1-(2-\mathbb{E}S) \frac{b_{N+1}-b_{N}}{a_{N}-a_{N+1}}
\end{equation*}
with large enough $N \in \mathbb{N}$. In all the examples below, we take
$N=20$. Using the same arguments as in Remark 2.1 of \cite{ds} we
can obtain both lower and upper bounds for $\psi (0)$ by calculating
$\psi _{N}(0)$ and $\psi _{N+1}(0)$. Then the upper bound for the
approximation error of $\psi (0)$ can be calculated by
\begin{equation*}
\Delta = |\psi _{N}(0)-\psi _{N+1}(0)|.
\end{equation*}

Finally, we can obtain approximations of the ruin probabilities\index{ruin probability} using
formula \eqref{a2} from Theorem~\ref{th1}
\begin{equation*}
1 - \psi (u) = a_{u} \bigl(1-\psi _{N}(0)\bigr) +
b_{u} (2-\mathbb{E}S), \quad u \in \mathbb{N}.
\end{equation*}

\begin{example}
\label{1ex}
Assume that the joint probability mass function of $(X,Y)$ is given by
the bivariate Poisson distribution:\index{bivariate Poisson distribution}
\begin{equation*}
\mathbb{P}(X=k, Y=l)=\sum_{i=0}^{\min \{k,l\}}
\frac{(\lambda _{1}-
\lambda )^{k-i}(\lambda _{2}-\lambda )^{l-i}\lambda ^{i}}{(k-i)!(l-i)!i!}
e ^{-(\lambda _{1}+\lambda _{2}-\lambda )}, \quad  k, l \in \mathbb{N}_{0},
\end{equation*}
where $\lambda _{j}>0$, $j=1,2$, $0\leqslant \lambda <\min \{\lambda
_{1},\lambda _{2}\}$. Then the marginal distribution\index{marginal distribution} of $X$ is Poisson\index{Poisson}
with parameter $\lambda _{1}$, the marginal distribution\index{marginal distribution} of $Y$ is
Poisson\index{Poisson} with parameter $\lambda _{2}$, and $\mathrm{{Cov}}(X,Y)=\lambda
$. If $\lambda =0$, then the two variables are independent, and the
results in this case are obtained in \cite{ds}.

In this example, we take $\lambda _{1}=0.3$ and $\lambda _{2}=1.4$. We
consider three possible values for the covariance parameter
$\lambda = \{0.01; 0.15; 0.29\}$, and the corresponding correlations
equal $\{0; 0.23; 0.46\}$.

In the table and graph below, the results of simulation are given. The
ruin probability\index{ruin probability} is calculated for the three values of the covariance
parameter mentioned above, and the upper bounds for the approximation
errors of $\psi (0)$ are also given.

From the results of simulation it could be observed, that for positively
dependent claims the ruin probability\index{ruin probability} is decreasing more slowly. It is
also interesting to note that the value of $\psi (0)$ is largest in the
case of independent claims.

\begin{table}[t]
\label{1table}
%
\caption{Values of $\psi (u)$ in Example~\ref{1ex}}
\begin{tabular}{rrrr}
\hline
$u$ & $\mathrm{{cor}}=0$ $(\Delta < 10^{-11})$
& $\mathrm{cor}=0.23$ $(\Delta < 10^{-10})$
& $\mathrm{cor}=0.46$ $(\Delta < 10^{-9})$ \\
\hline
$0$ & $0.7977$ & $0.7921$ & $0.7868$ \\
$1$ & $0.6040$ & $0.6264$ & $0.6480$ \\
$2$ & $0.4469$ & $0.4875$ & $0.5222$ \\
$3$ & $0.3269$ & $0.3754$ & $0.4165$ \\
$4$ & $0.2383$ & $0.2880$ & $0.3310$ \\
$5$ & $0.1736$ & $0.2208$ & $0.2628$ \\
$6$ & $0.1265$ & $0.1692$ & $0.2085$ \\
$7$ & $0.0921$ & $0.1297$ & $0.1655$ \\
$8$ & $0.0671$ & $0.0994$ & $0.1313$ \\
$9$ & $0.0489$ & $0.0762$ & $0.1042$ \\
$10$ & $0.0356$ & $0.0584$ & $0.0827$ \\
$11$ & $0.0260$ & $0.0447$ & $0.0657$ \\
$12$ & $0.0189$ & $0.0343$ & $0.0521$ \\
\hline
\end{tabular}
%
%
\vspace*{12pt}
\end{table}

\begin{figure}
\includegraphics{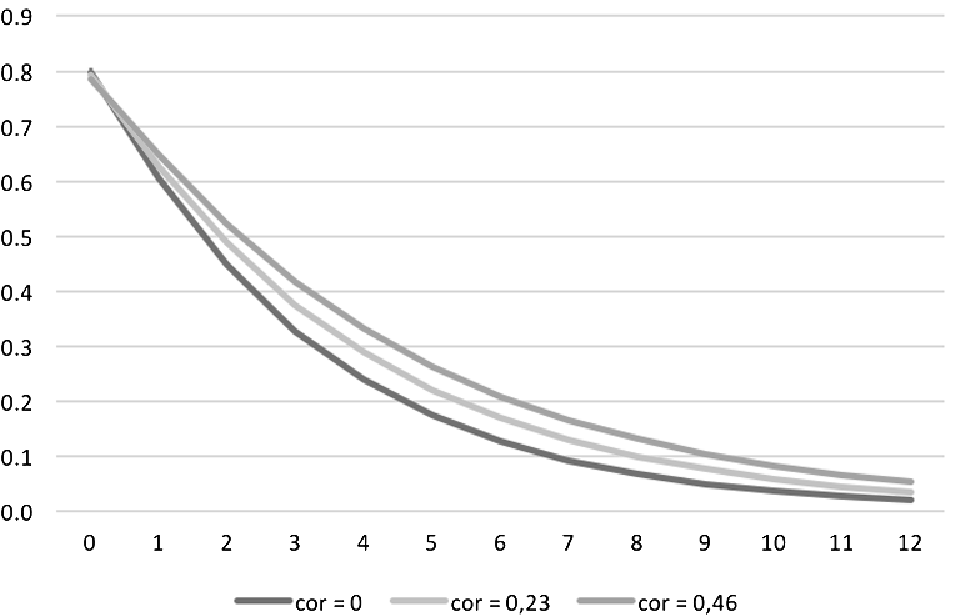}
\caption{Values of $\psi (u)$ in Example~\ref{1ex}}
\end{figure}
\end{example}

\begin{example}
\label{2ex}
This example deals with a Clayton copula\index{Clayton copula} and Poisson marginals.\index{Poisson marginals} Let us
denote $u_{1} := F_{X}(x)$, $u_{2} := F_{Y}(y)$. Clayton copula\index{Clayton copula} is
defined by
\begin{equation*}
C(u_{1},u_{2};\theta ) = \max \bigl\{u_{1}^{-\theta }
+ u_{2}^{-\theta } - 1, 0\bigr\}^{-{1}\backslash {\theta }}, \quad
u_{1}, u_{2} \in [0, 1],
\end{equation*}
where the dependence parameter $\theta \in [-1, \infty ) \backslash
\{0\}$. The marginals become independent as $\theta \rightarrow 0$.
Clayton copula\index{Clayton copula} can be used to model negative dependence\index{negative dependence} when
$\theta \in [-1, 0)$. Detailed analysis of this copula can be found, for
instance, in \cite{j-2014,ml,mn-2009} and \cite{n-2006}.

In this example, the marginal distribution\index{marginal distribution} of $X$ is Poisson\index{Poisson} with
parameter $0.3$, and the marginal distribution\index{marginal distribution} of $Y$ is Poisson\index{Poisson} with
parameter $1.4$. We take three values for the covariance parameter
$\theta = \{-0.9; 0.01; 100\}$, and the corresponding correlations equal
$\{-0.53; 0; 0.8\}$.\vadjust{\eject}

From the results of simulation it could be observed, that as in Example~\ref{1ex}
 for positively dependent claims the ruin probability\index{ruin probability} is
decreasing more slowly. It is also interesting to note that the value
of $\psi (0)$ is largest in the case of negatively dependent claims.

\begin{table}[t]
\label{2table}
%
\caption{Values of $\psi (u)$ in Example~\ref{2ex}}
\begin{tabular}{rrrr}
\hline
$u$ & $\mathrm{cor}=-0.53$ $(\Delta < 10^{-20})$
& $\mathrm{cor}=0$ $(\Delta < 10^{-11})$
& $\mathrm{cor}=0.8$ $(\Delta < 10^{-10})$ \\
\hline
$0$ & $0.8217$ & $0.7977$ & $0.7810$ \\
$1$ & $0.5064$ & $0.6040$ & $0.6717$ \\
$2$ & $0.3165$ & $0.4469$ & $0.5715$ \\
$3$ & $0.1977$ & $0.3269$ & $0.4669$ \\
$4$ & $0.1231$ & $0.2383$ & $0.3909$ \\
$5$ & $0.0766$ & $0.1736$ & $0.3221$ \\
$6$ & $0.0476$ & $0.1265$ & $0.2661$ \\
$7$ & $0.0296$ & $0.0921$ & $0.2195$ \\
$8$ & $0.0184$ & $0.0671$ & $0.1812$ \\
$9$ & $0.0115$ & $0.0489$ & $0.1496$ \\
$10$ & $0.0071$ & $0.0356$ & $0.1235$ \\
$11$ & $0.0044$ & $0.0260$ & $0.1019$ \\
$12$ & $0.0028$ & $0.0189$ & $0.0841$ \\
\hline
\end{tabular}
%
%
\end{table}

\begin{figure}
\includegraphics{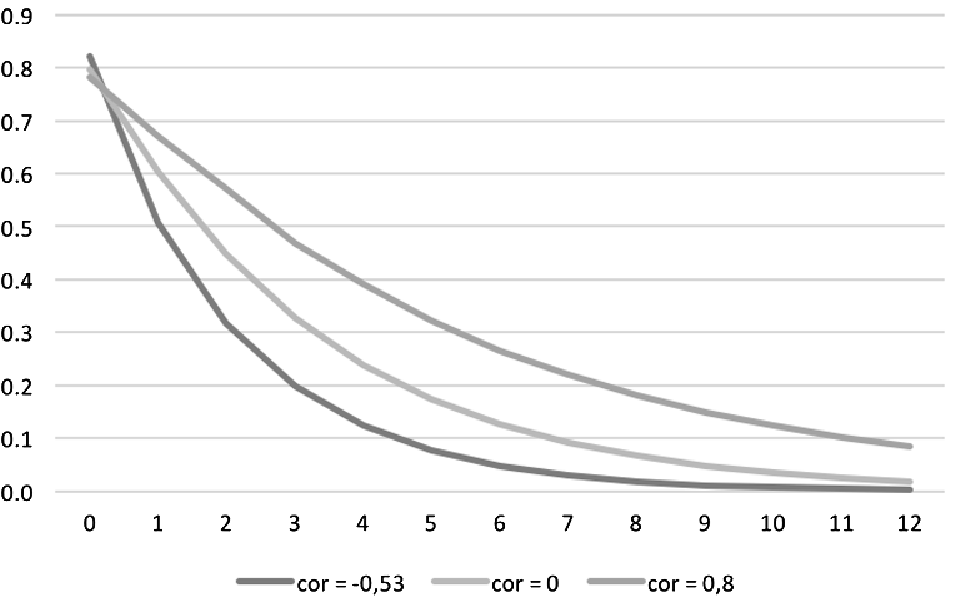}
\caption{Values of $\psi (u)$ in Example~\ref{2ex}}
\end{figure}
\end{example}

\begin{example}
\label{3ex}
This example is the opposite case of Example~\ref{2ex}. The marginal
distribution\index{marginal distribution} of $X$ is Poisson\index{Poisson} with parameter $1.4$, and the marginal
distribution\index{marginal distribution} of $Y$ is Poisson\index{Poisson} with parameter $0.3$. To model the
dependence between the marginals,\index{marginal} we use the Clayton copula\index{Clayton copula} with
$\theta = \{-0.9; 0.01; 100\}$ again, and the corresponding correlations
equal $\{-0.53; 0; 0.8\}$.

From the simulation we can observe that the order of appearance of
claims has considerable effect on the ruin probability.\index{ruin probability}

\begin{table}[t]
\label{3table}
%
\caption{Values of $\psi (u)$ in Example~\ref{3ex}}
\begin{tabular}{rrrr}
\hline
$u$ & $\mathrm{cor}=-0.53$ $(\Delta < 10^{-20})$
& $\mathrm{cor}=0$ $(\Delta < 10^{-11})$
& $\mathrm{cor}=0.8$ $(\Delta < 10^{-9})$ \\
\hline
$0$ & $0.9267$ & $0.9023$ & $0.8988$ \\
$1$ & $0.6940$ & $0.7269$ & $0.7316$ \\
$2$ & $0.4653$ & $0.5473$ & $0.5897$ \\
$3$ & $0.2961$ & $0.4014$ & $0.4859$ \\
$4$ & $0.1850$ & $0.2926$ & $0.4048$ \\
$5$ & $0.1151$ & $0.2131$ & $0.3347$ \\
$6$ & $0.0716$ & $0.1552$ & $0.2763$ \\
$7$ & $0.0445$ & $0.1131$ & $0.2280$ \\
$8$ & $0.0277$ & $0.0824$ & $0.1882$ \\
$9$ & $0.0172$ & $0.0600$ & $0.1553$ \\
$10$ & $0.0107$ & $0.0437$ & $0.1282$ \\
$11$ & $0.0067$ & $0.0319$ & $0.1059$ \\
$12$ & $0.0042$ & $0.0232$ & $0.0874$ \\
\hline
\end{tabular}
%
%
\end{table}

\begin{figure}
\includegraphics{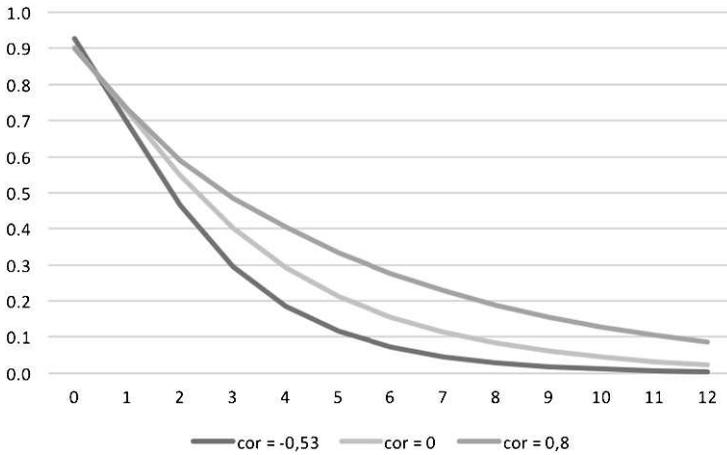}
\caption{Values of $\psi (u)$ in Example~\ref{3ex}}
\end{figure}
\end{example}

\begin{example}
\label{4ex}
All the examples considered so far deal only with light-tailed
marginals, but Theorem~\ref{th1} only imposes requirement for the
expectations of the marginals\index{marginal} while higher order moments can be
infinite. In this example, the distribution of the first claim $X$ is
Poisson\index{Poisson} with parameter $\lambda =0.2$, and the second claim $Y$ is
distributed according to the Zeta distribution\index{Zeta distribution} with parameter $2.3$,
that is
\begin{equation*}
\mathbb{P}(Y=m)=\frac{1}{\zeta (2.3)}\frac{1}{(m+1)^{2.3}},\quad  m\in \mathbb{N}_{0},
\end{equation*}
where $\zeta $ denotes the Riemann zeta function. It should be noted
that here Zeta distribution\index{Zeta distribution} is not defined in the usual way, i.e. with
support $m \in \{1,2,\ldots \}$ and the corresponding probabilities.

The expectation of $Y$ is $1.74497$ and the variance is infinite.
Therefore, the correlation between the claims is undefined. As before,
we use the Clayton copula\index{Clayton copula} with $\theta = \{-0.9; 0.01; 100\}$ to model
the dependence between the marginals.\index{marginal}

As can be intuitively expected, the presence of heavy-tailed marginal
has a major impact on the values of the ruin probability.\index{ruin probability}

\begin{table}[t]
\label{4table}
%
\caption{Values of $\psi (u)$ in Example~\ref{4ex}}
\begin{tabular}{rrrr}
\hline
$u$ & $\theta =-0.9$ $(\Delta < 10^{-6})$
& $\theta =0.01$ $(\Delta < 10^{-6})$
& $\theta =100$ $(\Delta < 10^{-5})$ \\
\hline
$0$ & $0.9721$ & $0.9715$ & $0.9690$ \\
$1$ & $0.9611$ & $0.9620$ & $0.9656$ \\
$2$ & $0.9570$ & $0.9579$ & $0.9615$ \\
$3$ & $0.9543$ & $0.9550$ & $0.9584$ \\
$4$ & $0.9520$ & $0.9527$ & $0.9559$ \\
$5$ & $0.9500$ & $0.9507$ & $0.9538$ \\
$6$ & $0.9483$ & $0.9489$ & $0.9520$ \\
$7$ & $0.9467$ & $0.9473$ & $0.9503$ \\
$8$ & $0.9453$ & $0.9458$ & $0.9488$ \\
$9$ & $0.9439$ & $0.9444$ & $0.9474$ \\
$10$ & $0.9427$ & $0.9432$ & $0.9460$ \\
$11$ & $0.9416$ & $0.9421$ & $0.9448$ \\
$12$ & $0.9406$ & $0.9410$ & $0.9437$ \\
\hline
\end{tabular}
%
%
\end{table}

\begin{figure}
\includegraphics{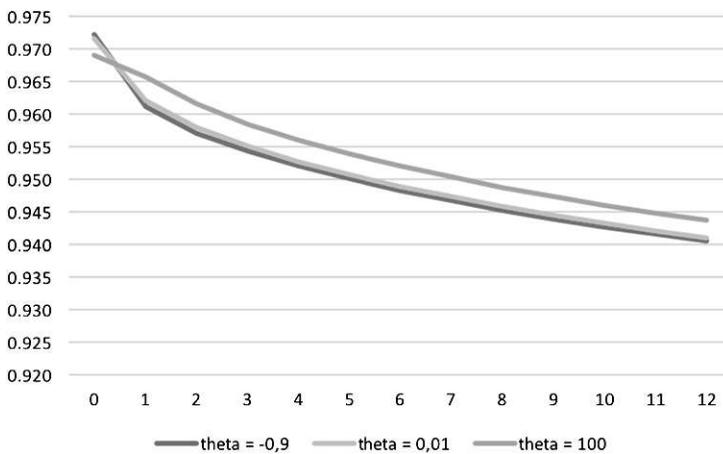}
\caption{Values of $\psi (u)$ in Example~\ref{4ex}}
\end{figure}
\end{example}

\section{Concluding remarks}%
\label{concl}
In this work, the bi-seasonal discrete time risk model\index{discrete time risk model} with dependent
claims is introduced. We present a recursive algorithm for calculating
the values of the ruin probability.\index{ruin probability} Theoretical results are illustrated
by some numerical examples.

The results obtained in this paper can be extended in the following
directions:
\begin{itemize}
\item
Our results can be generalized to the models with more complex structure
of the non-homogeneity of claims. For instance, the generating random
vectors of the form $(X_{1}, X_{2}, \ldots , X_{p})$ with $p>2$ can be
considered for claim sizes. In this case, we get a $p$-seasonal model.
\item
An algorithm for the calculation of more complex risk measures, such as
the Gerber--Shiu expected discounted penalty function
\cite{gs-1998}, can be presented for the bi-seasonal discrete time risk
model\index{discrete time risk model} with dependent claims.
\item
The model and the algorithm considered in the paper can be illustrated
with examples based on real insurance data.
\end{itemize}


\begin{acknowledgement}
We are grateful to the referees for their
useful comments and suggestions leading to an improvement of the paper.
\end{acknowledgement}

\begin{funding}
The second  and the third authors were supported by grant No \gnumber[refid=GS1]{S-MIP-17-72} from
the \gsponsor[id=GS1,sponsor-id=501100004504]{Research
Council of Lithuania}.
\end{funding}


\end{document}